\documentstyle[12pt]{article}
\def\Bbb R{{\rm \bf R}}
\def\proclaim#1{\vskip2mm{\bf #1}\em}
\def\endproclaim{\em \vskip2mm}
\def\tag#1{\eqno(#1)}
\def\gathered{\begin{array}{c}}
\def\endgathered{\end{array}}
\def\text{\mbox}

\begin{document}

\title {The enclosure method for inverse obstacle scattering over a finite time interval: IV.
Extraction from a single point on the graph of the response operator}
\author{Masaru IKEHATA\footnote{
Laboratory of Mathematics,
Institute of Engineering,
Hiroshima University,
Higashi-Hiroshima 739-8527, JAPAN}}
\maketitle

\begin{abstract}
Now a final and maybe simplest formulation of the enclosure method
applied to inverse obstacle problems governed by partial differential equations
in a {\it spacial domain with an outer boundary} over a finite time interval
is fixed.  The method employs only a single pair of a quite natural Neumann data
prescribed on the outer boundary and the corresponding Dirichlet data on the same boundary
of the solution of the governing equation
over a finite time interval, that is a single point on the graph of the so-called {\it response operator}.
It is shown that the methods enables us to extract the distance of a given point outside the domain to an embedded unknown obstacle,
that is the maximum sphere centered at the point whose exterior encloses the unknown obstacle.
To make the explanation of the idea clear only an inverse obstacle problem governed by the wave equation is considered.

\noindent
AMS: 35R30, 35L05

\noindent KEY WORDS: enclosure method, inverse obstacle problem, wave equation, heat equation,
non destructive testing.
\end{abstract}


\section{Introduction}

As done in the series of the previous papers \cite{IW00, IEO2, IEO3} the aim of this paper is to pursuit the possibility of the enclosure method itself for inverse obstacle problems in time domain.  
This paper adds an {\it extremly simple} method employing  the enclosure method as a guiding principle
to the list of previous versions of the enclosure method.
It is rigorous and applicable to a broad class of inverse obstacle problems governed by partial differential equations in time domain, 
including heat and wave equations in a domain {\it with an outer boundary}
over a finite time interval.
Such class should be a mathematical counterpart of, fo example, a non destructive testing using acoustic and elastic waves in time domain.

Now let us descibe the simple method mentioned above.
To show the idea clearly we restrict ourself to
an inverse obstacle problem using a scalar wave which propagates inside a three dimensional body.

Let $\Omega$ be a bounded domain with $C^2$-boundary.
Le $D$ be a nonempty bounded open subset of $\Omega$ with $C^2$-boundary such that
$\Omega\setminus\overline D$ is connected.
Let $0<T<\infty$.

Given $f=f(x,t), (x,t)\in\partial\Omega\times\,]0,\,T[$ let $u=u_f(x,t), (x,t)\in\,(\Omega\setminus\overline D)\times\,]0,\,T[$ 
denote the solution of the following initial boundary value problem for the classical wave equation:
$$\displaystyle
\left\{
\begin{array}{ll}
\displaystyle
(\partial_t^2-\Delta) u=0 & \text{in}\,(\Omega\setminus\overline D)\times\,]0,\,T[,\\
\\
\displaystyle
u(x,0)=0 & \text{in}\,\Omega\setminus\overline D,
\\
\\
\displaystyle
\partial_tu(x,0)=0 & \text{in}\,\Omega\setminus\overline D,\\
\\
\displaystyle
\frac{\partial u}{\partial\nu}=0 & \text{on}\,\partial D\times\,]0,\,T[,\\
\\
\displaystyle
\frac{\partial u}{\partial\nu}=f(x,t) & \text{on}\,\partial\Omega\times\,]0,\,T[.
\end{array}
\right.
\tag {1.1}
$$
We use the same symbol $\nu$ to denote both the outer unit normal vectors of $\partial D$
and $\partial\Omega$. 
The solution class and the Neumann data $f$ should be specified later.

We consider the following problem.

{\bf\noindent Problem.}  Fix a large $T$ (to be determined later) and a single $f$ (to be specified later).
Assume that set $D$ is  unknown.  Extract information about the
location and shape of $D$ from the wave field $u_f(x,t)$ given at
all $x\in\partial\Omega$ and $t\in\,]0,\,T[$.

As called in the BC-method \cite{B}, the correspondence $f\longmapsto u_f\vert_{\partial\Omega\times\,]0,\,T[}$ 
should be called the {\it response operator}.  However, unlike the BC-method, we try to extract some information
about the geometry of unknown obstacle from $u_f$ on $\partial\Omega\times\,]0,T[$ for a {\it fixed} $f$,
that is a point on the graph of the response operator.

Since (1.1) is a {\it non-homogeneous} Neumann problem, the solution class for {\it general} Neumann data
$f$ is not simple compared with the {homogeneous Neumann problem} which can be covered by
a variational approach \cite{DL} or the theory of $C_0$-semigroup \cite{Y} in the $L^2$-frame work.
See \cite{OK2}, therein a {\it fractional} Sobolev space is used for the description of the solution
class for the nonhomogeneous Neumann problem for applying the BC-method.
Then, in this paper, we do not prescribe the completely general Neumann data $f$, instead, generate the necessary $f$
by solving the wave equation in the whole space.

Let $B$ be an open ball satisfying $\overline B\cap\overline\Omega=\emptyset$.
We think the radius $\eta$ of $B$ is very small.
Let $\chi_B$ denote the characteristic function of $B$.  Let $v=v_B$ solve
$$\displaystyle
\left\{
\begin{array}{ll}
\displaystyle
(\partial_t^2-\Delta)v=0 & \text{in}\,\Bbb R^3\times\,]0,\,T[,\\
\\
\displaystyle
v(x,0)=0 & \text{in}\,\Bbb R^3,\\
\\
\displaystyle
\partial_tv(x,0)=\Psi_B(x) & \text{in}\,\Bbb R^3,
\end{array}
\right.
\tag {1.2}
$$
where 
$$\begin{array}{ll}
\displaystyle
\Psi_B(x)=(\eta-\vert x-p\vert)\chi_B(x), & x\in\Bbb R^3
\end{array}
\tag {1.3}
$$
and $p$ denotes the center of $B$.
Note that the function $\Psi_B$ belongs to $H^1(\Bbb R^3)$ since
$\nabla\Psi_B(x)=-\{(x-p)/\vert x-p\vert\}\chi_B(x)$ in the sense of distribution.
The solution $v_B$ of (1.2) is constructed by using the theory of $C_0$-semigroupe.
The class where $v_B$ belongs to is the following:
$$\displaystyle
v_B\in C^2([0,\,T], L^2(\Bbb R^3))\cap C^1([0,\,T], H^1(\Bbb R^3))\cap C([0,\,T], H^2(\Bbb R^3)).
$$
Needless to say, $v_B$ has an explicit analytical expression, however, we never make use of 
such expression in time domain.  We need just the exsitence of $v_B$ in the function spaces indicated above.

The following function is the special $f$ in the problem mentioned above.

Define
$$\begin{array}{ll}
\displaystyle
f_B=f_B(\,\cdot\,,t)=\frac{\partial}{\partial\nu}v_B(\,\cdot\,,t), & t\in\,[0,\,T].
\end{array}
\tag {1.4}
$$
Note that function $f_B$ does not contain any {\it large} parameter.

Now we construct the solution of (1.1) by prescribing $f=f_B$.
First we make use of a standard reduction of non-homogeneous Neumann problem
to homogeneous one by using the special form (1.4) of the Neumann data.

Since $\partial\Omega$ is $C^2$, one can choose a $C^2$-function $\phi$ such that 
$\phi=1$ in a neighbourhood of $\partial\Omega$ and $\phi=0$ in a neighbourhood of $\overline D$ and
the outside of an open ball with a large radius containing $\overline\Omega$.
We have
$$\displaystyle
(\partial_t^2-\Delta)(\phi v_B)=-(\Delta\phi)v_B-2\nabla\phi\cdot\nabla v_B\in C^1([0,T], L^2(\Bbb R^3)).
$$
Then, by applying the theory of $C_0$-semigroup, we have the unique
$z\in C^2([0, T], L^2(\Omega\setminus\overline D))\cap C^1([0,\,T], H^1(\Omega\setminus\overline D))\cap C([0,\,T], H^2(\Omega\setminus\overline D))$
such that
$$\displaystyle
\left\{
\begin{array}{ll}
\displaystyle
(\partial_t^2-\Delta)z=(\partial_t^2-\Delta)(\phi v_B) & \text{in}\,(\Omega\setminus\overline D)\times\,]0,\,T[,\\
\\
\displaystyle
z(x,0)=0 & \text{in}\,\Omega\setminus\overline D,
\\
\\
\displaystyle
\partial_tz(x,0)=0 & \text{in}\,\Omega\setminus\overline D,\\
\\
\displaystyle
\frac{\partial z}{\partial\nu}=0 & \text{on}\,\partial D\times\,]0,\,T[,\\
\\
\displaystyle
\frac{\partial z}{\partial\nu}=0 & \text{on}\,\partial\Omega\times\,]0,\,T[.
\end{array}
\right.
$$
We refer the reader to Theorem 1 in \cite{I} which includes more general homogeneous boundary condition.
Then, the $u$ defined by
$$\displaystyle
u=\phi v_B-z\in C^2([0, T], L^2(\Omega\setminus\overline D))\cap C^1([0,\,T], H^1(\Omega\setminus\overline D))
\cap C([0,\,T], H^2(\Omega\setminus\overline D)),
$$
is the desired solution of (1.1).  The uniquness in this class is clear.

Now having the solution $u=u_f$ of (1.1) with $f=f_B$ given by (1.4), we set
$$\begin{array}{lll}
\displaystyle
w_B(x)=w_B(x,\tau)=\int_0^Te^{-\tau t}u_f(x,t)dt,
&
x\in\Omega\setminus\overline D,
&
\tau>0
\end{array}
\tag {1.5}
$$
and
$$\begin{array}{lll}
\displaystyle
w_B^0(x)=w_B^0(x,\tau)=\int_0^Te^{-\tau t}v_B(x,t)dt,
&
x\in\Bbb R^3,
&
\tau>0.
\end{array}
\tag {1.6}
$$

Define
$$
\begin{array}{ll}
\displaystyle
I_{\partial\Omega}(\tau;B)=\int_{\partial\Omega}
(w_B-w_B^0)\frac{\partial w_B^0}{\partial\nu}\,dS, & \tau>0.
\end{array}
\tag {1.7}
$$
This is the indicator function in the {\it enclosure method} developed in this paper.

This indicator function can be computed from the responce $u_B$ on $\partial\Omega$ over
time interval $]0,\,T[$ which is the solution of (1.1) with $f=f_B$.

Now we state the main result of this paper.

\proclaim{\noindent Theorem 1.1.}

(i)  If $T$ satisfies
$$\displaystyle
T>2\text{dist}\,(D,B)-\text{dist}\,(\Omega,B),
\tag {1.8}
$$
then, there exists a positive number $\tau_0$ such that, for all $\tau\ge\tau_0$
$I_{\partial D}(\tau;B)>0$ and we have
$$\displaystyle
\lim_{\tau\longrightarrow\infty}
\frac{1}{\tau}
\log I_{\partial\Omega}(\tau;B)=-2\text{dist}\,(D,B).
\tag {1.9}
$$

(ii)  We have
$$\displaystyle
\lim_{\tau\longrightarrow\infty}e^{\tau T}I_{\partial\Omega}(\tau;B)=
\left\{
\begin{array}{ll}
\displaystyle
\infty & \text{if}\,\,T>2\text{dist}\,(D,B),\\
\\
\displaystyle
0     & \text{if}\,\,T<2\text{dist}\,(D,B).
\end{array}
\right.
\tag {1.10}
$$

\endproclaim

Note that if $T=2\text{dist}\,(D,B)$, the proof tells us only $e^{\tau T}I_{\partial\Omega}(\tau;B)=O(\tau^4)$ 
as $\tau\longrightarrow\infty$.  And also note that we omitted to denote the dependence of $I_{\partial\Omega}(\tau;B)$
on $T$.

In short, Theorem 1.1 says that the out put generated by a single input depending on $B$ and given on the boundary 
of the domain over a finite time interval uniquely determines $\text{dist}\,(D,B)$.
Define $d_{\partial D}(p)=\inf_{y\in\partial D}\vert y-p\vert$.  We have
$\text{dist}\,(D,B)=d_{\partial D}(p)-\eta$.  Since $B$ is known, we can conclude that the indicator function for each $B$
uniquely determines $d_{\partial D}(p)$ and hence the sphere $\vert x-p\vert=d_{\partial D}(p)$ on which there exists a point
on $\partial D$.  This sphere is the maxium one whose exterior contains $D$.   Moving $p$ outside $\Omega$, we can obtain
an estimation of the geometry of $D$.  The point is, one input yields one information.  We do not need the whole knowledge
of the response operator before doing the procedure.

The restriction (1.8) is an effect on the measurement on $\partial\Omega$.  Note also that as pointed out in \cite{IW00}
we have the inequality:
$$\displaystyle
2\text{dist}\,(D,B)-\text{dist}\,(\Omega,B)\ge
\inf\left\{\vert x-y\vert+\vert y-z\vert\,\vert\,x\in\partial B, y\in\partial D,\,z\in\partial\Omega
\right\}.
$$
From a geometrical optics point of view the quantity on this right-hand side can be interpreted as the first arrival time of
a {\it virtual signal} that strarts from the surface of $B$ at $t=0$, reflects on the surface of the obstacle
and arrives at a point on $\partial\Omega$.  Note that the solution of (1.1) describes a wave which propagates
the spatial domain $\Omega\setminus\overline D$ only.
However, as can be seen in the definition of the indicator function,
we generate a wave inside domain $\Omega\setminus\overline D$ by using the special Neumann data $f_B$ on $\partial\Omega$ over finite time interval
$]0,\,T[$ given by (1.3).  Theorem 1.1 suggests us the design of the Neumann data that makes the boundary
of domain $\Omega$ {\it transparent} and enables us to extract the distance of $D$
from $B$ directly.  Note that the obtained quantity $\text{dist}\,(D,B)$ is simpler
than the quantity mentioned above and enables us easily to find an estimation of the location of unknown obstacle from above. 
In practice, we should develop a realization method of the Neumann data desired in Theorem 1.1
by using a principle of superposition.

It follows from (ii) in Theorem 1.1 that the formula
$$\displaystyle
2\text{dist}\,(D,B)
=\sup\left\{T\in\,]0,\,\infty[\,\vert\,\lim_{\tau\longrightarrow\infty}e^{\tau T}I_{\partial\Omega}(\tau;B)=0
\right\},
$$
is valid.  This formula has a similarity with an original version of the enclosure method in \cite{E00}.
See also (1.13) in \cite{IMax2} fot the Maxwell system in an exterior domain.

Some corollaries of Theorem 1.1 are in order.
Fist let $v_0\in H^1(\Bbb R^3)$ be the weak solution of 
$$\begin{array}{ll}
\displaystyle
(\Delta-\tau^2)v+\Psi_B=0 & \text{in}\,\Bbb R^3.
\end{array}
\tag {1.11}
$$
Define another indicator function
$$
\begin{array}{ll}
\displaystyle
I_{\partial\Omega}^s(\tau;B)=\int_{\partial\Omega}
(w_B-w_B^0)\frac{\partial v_0}{\partial\nu}\,dS, & \tau>0.
\end{array}
\tag{1.12}
$$

Note that $v_0$ has the expression
$$\displaystyle
v_0(x)=\frac{1}{4\pi}\int_B\frac{e^{-\tau\vert x-y\vert}}{\vert x-y\vert}\cdot(\eta-\vert y-p\vert)dy.
\tag {1.13}
$$
Thus indicator function $I^s_{\partial\Omega}(\tau;B)$ is simpler than the former indicator function $I_{\partial\Omega}(\tau;B)$.

\proclaim{\noindent Corollary 1.1.}
All the statements in Theorem 1.1 for $I_{\partial\Omega}(\tau;B)$ 
replaced with $I_{\partial\Omega}^s(\tau;B)$, are valid.

\endproclaim

Finally we introduce a {\it localization} of indicator function $I_{\partial\Omega}^s(\tau;B)$.
Given $M>0$ define
$$\displaystyle
\partial\Omega(B,M)=\left\{x\in\partial\Omega\,\vert\,d_B(x)<M\right\},
$$
where $d_B(x)=\inf_{y\in B}\vert y-x\vert$.

Define the {\it localized indicator function} $I_{\partial\Omega}(\tau;B,M)$ by the formula
$$\begin{array}{ll}
\displaystyle
I_{\partial\Omega}(\tau;B,M)
=\int_{\partial\Omega(B,M)}
(w_B-w_B^0)\frac{\partial v_0}{\partial\nu}\,dS, & \tau>0.
\end{array}
$$
We are ready to state the second corollay of Theorem 1.1.

\proclaim{\noindent Corollary 1.2.}
Let $M$ satisfy
$$\displaystyle
\text{dist}\,(D,B)<M.
\tag {1.14}
$$
Let $T$ satisfy
$$\displaystyle
T\ge 2M-\text{dist}\,(\Omega,B).
\tag {1.15}
$$
Then, the statement (i) in Theorem 1.1 for $I_{\partial\Omega}(\tau;B)$ 
replaced with $I_{\partial\Omega}(\tau;B,M)$, is valid.

\endproclaim

The $M$ in (1.14) plays a role of a-priori information about the location of $D$ from $B$.
Corollary 1.2 shows that with the help of this information one can reduce the size of the place where the data are collected.

\subsection{Comparison with the previous enclosure method in time domain}

The enclosure method for inverse obstacle problems in time domain was initiated in \cite{I4}
and its idea goes back to the method developed in \cite{I1}.
In \cite{I4} the author considered some prototype inverse obstacle problems for the heat equation in one-space dimensional case
and found the enclosure method using a {\it single set} of lateral data over a finite time interval.
The method makes use of a special solution of a formal adjoint of the governing equation for the background medium 
or related equation depending on a large 
parameter often denoted by $\tau$ and observation data. Using integration by parts, from those we construct
an {\it indicator function} of indepenent variable $\tau$.
From the asymptotic behaviour of the indicator function as $\tau\longrightarrow\infty$ we find a domain that encloses unknown obstacles.
This idea is realized in three-space dimensions for inverse obstacle problems governed by the wave equations \cite{IW00, ICA,
IEO2, IEO3, IEE, Iwall}, the Maxwell system \cite{IMax, IMax2} and heat equations \cite{IK2, IK5}.

It is worth comparing the method in this paper with the methods in \cite{IW00, IEO2}.
One of the inverse obstacle problems considered therein is the following.
Consider the following initial exterior boundary value problem:
$$\displaystyle
\left\{
\begin{array}{ll}
\displaystyle
(\partial_t^2-\Delta) u=0 & \text{in}\,(\Bbb R^3\setminus\overline D)\times\,]0,\,T[,\\
\\
\displaystyle
u(x,0)=0 & \text{in}\,\Bbb R^3\setminus\overline D,\\
\\
\displaystyle
\partial_tu(x,0)=\chi_B(x) & \text{in}\,\Bbb R^3\setminus\overline D,\\
\\
\displaystyle
\frac{\partial u}{\partial\nu}=0 & \text{on}\,\partial D\times\,]0,\,T[.
\end{array}
\right.
\tag {1.16}
$$
Note that $B$ is an open ball satisfying $\overline B\cap\overline D=\emptyset$.
As the measurement place we choose a bounfded open set $\Omega'$ of $\Bbb R^3$ with a smooth boundary
satisfying $\overline B\cap\overline{\Omega'}=\emptyset$ and $\overline D\subset\Omega'$.
We denote by $\nu$ again the unit outward normals to $\partial D$ and $\partial\Omega'$.
The inverse problem is to extract information about the location and shape of
$D$ from $u$ and $\partial u/\partial\nu$ on $\partial\Omega'\times\,]0\,T[$ for a fixed large $T$ and $B$.
The method developed therein is the following.

Let $\tau>0$.  Let $v_0'\in H^1(\Bbb R^3)$ be the weak solution of
$$\begin{array}{ll}
\displaystyle
(\Delta-\tau^2)v+\chi_B=0 & \text{in $\Bbb R^3$.}
\end{array}
\tag {1.17}
$$
We introduced the indicator function by the formula
$$
\displaystyle
I'_{\partial\Omega'}(\tau;B)
=\int_{\partial\Omega'}
\left(w'\frac{\partial v_0'}{\partial\nu}-v_0'\frac{\partial w'}{\partial\nu}\right)\,dS,
$$
where
$$\begin{array}{ll}
\displaystyle
w'=w'(x,\tau)
=\int_0^T e^{-\tau t}u'(x,t)dt, & x\in\Bbb R^3\setminus\overline D,
\end{array}
$$
and $u'$ is the solution of (1.16).  Under the assumption (1.8) in which $\Omega$ is replaced with $\Omega'$ we obtained
a formula corresponding to formula (1.9).  In addtion, as can be seen in a recent application \cite{Iwall} of the enclosure method for 
inverse obstacle problems arising in {\it through-wall imaging} one can replace $v_0'$ with the function
$$\displaystyle
\int_0^{T}e^{-\tau t}v_B'(x,t)dt,
$$
where $v'=v_B'$ solves (1.2) with $\Psi_B$ replaced with $\chi_B$.  Thus a choice or generating method of a special
solution needed has a common point in the spirit. 
Note also that in \cite{IEO2}, we have pointed out that, as $\tau\longrightarrow\infty$
$$\displaystyle
I'_{\partial\Omega'}(\tau;B)
=\int_B(w'-v_0')dx+O(\tau^{-1}e^{-\tau T}).
$$
Using this relationship we have transplanted all the results for $I'_{\partial\Omega'}(\tau;B)$
into those for another indicator function
$$\displaystyle
\tau\longmapsto\int_B(w'-v_0')dx,
$$
which can be computed by using the back-scattering data $u'(x,t)$ given at all $x\in B$ and $t\in\,]0,\,T[$.

From the comparison above, in short, in this paper we have found a counterpart of the methods developed in \cite{IW00, IEO2}
in a class of the inverse obstacle problems in time domain goverened by
partial differential equations defined in a {\it spacial domain with an outer boundary}
over a finite time interval.

Note that another enclosure method originating from \cite{E00}
and using {\it infinitely many} sets of lateral 
data over a finite time interval has been developed in \cite{IK1, IK2, IFR}
for the heat equations in three-space dimensions and parabolic system \cite{II}.
See also \cite{OK, OK2} which are based on the BC-method \cite{B} using the full knowledge of the response operator itself
for the wave equation over a finite time interval.
However, in this paper we employ only a singe point on the graph of the response operator
and so we will not discuss those methods here.

Finally we compare the method in this paper with a result in \cite{IK2} for the heat equation.  
Therein we considered an inverse initial
boundary value problem for the heat equation $(\partial_t-\nabla\cdot\gamma\nabla)u=0$
in $\Omega\times\,]0,\,T[$ with discontiunus coefficient $\gamma$ and an arbitrary fixed $T>0$. 
One of the results is:  
there is a computation formula of the distance of an unknown inclusion and $\partial\Omega$ from a single set of the temparture 
generated by a single input heat flux $f$ across $\partial\Omega$ over time interval $]0,\,T[$
under, roughly speaking, {\it positivity} of $f$ on $\partial\Omega\times\,]0,\,T[$.

Instead of $v'_0$ satisfying (1.17), the approach in \cite{IK2} employs the solution
of the nonhomgeneous Neumann problem
$$\left\{
\begin{array}{ll}
\displaystyle
(\Delta-\tau)v=0 & \text{in $\Omega$},\\
\\
\displaystyle
\frac{\partial v}{\partial\nu}=g & \text{on $\partial\Omega$},
\end{array}
\right.
\tag {1.18}
$$
where
$$\displaystyle
g=g(x,\tau)=\int_0^Te^{-\tau t}f(x,t)dt.
$$
The proof is based on the asymptotic behaviour of
the solution of (1.18) as $\tau\longrightarrow\infty$ in a neighbourhood of the closure of the inclusion.
The analysis makes use of an expression of the solution constructed by solving an integral equation on $\partial\Omega$.
Note that, developing this approach for a parabolic sysytem has been left as an open problem, see \cite{II}.
However, it would be possible to apply the method presented in this paper to the problem and shall
be reported in forthcoming papers.

\section{Proof of Theorem 1.1 and Corollaries}

In this section, for simplicity of description we always write
$$\left\{
\begin{array}{l}
\displaystyle
w=w_B,\\
\\
\displaystyle
w_0=w_B^0,\\
\\
\displaystyle
R=w-w_0,
\end{array}
\right.
$$
where $w_B$ and $w_B^0$ are given by (1.5) and (1.6).

\subsection{A decomposition formula of the indicator function}

It follows from (1.1) that $w$ satisfies
$$\left\{
\begin{array}{ll}
\displaystyle
(\Delta-\tau^2)w=e^{-\tau T}F & \text{in}\,\Omega\setminus\overline D,\\
\\
\displaystyle
\frac{\partial w}{\partial\nu}=\frac{\partial w_0}{\partial\nu} & \text{on}\,\partial\Omega,\\
\\
\displaystyle
\frac{\partial w}{\partial\nu}=0 & \text{on}\,\partial D,
\end{array}
\right.
\tag {2.1}
$$
where
$$\begin{array}{ll}
\displaystyle
F=F(x,\tau)=\partial_tu(x,T)+\tau u(x,T), & x\in\Omega\setminus\overline D.
\end{array}
$$
Note that we have
$$\displaystyle
\Vert F\Vert_{L^2(\Omega\setminus\overline D)}=O(\tau).
\tag {2.2}
$$

It follows from (1.2) that the $w_0$ satisfies
$$
\begin{array}{ll}
\displaystyle
(\Delta-\tau^2)w_0+\Psi_B=e^{-\tau T}F_0 & \text{in}\,\Bbb R^3,
\end{array}
\tag {2.3}
$$
where
$$\begin{array}{ll}
\displaystyle
F_0=F_0(x,\tau)
=\partial_tv_B(x,T)+\tau v_B(x,T), & x\in\Bbb R^3.
\end{array}
$$
Note that we have
$$\displaystyle
\Vert F_0\Vert_{L^2(\Bbb R^3)}=O(\tau).
\tag {2.4}
$$

Then, integration by parts together with (2.1) and (2.3) yields
$$\begin{array}{l}
\,\,\,\,\,\,
\displaystyle
\int_{\partial\Omega}
\left(\frac{\partial w_0}{\partial\nu}w-\frac{\partial w}{\partial\nu}w_0\right)dS\\
\\
\displaystyle
=\int_{\partial D}w\frac{\partial w_0}{\partial\nu}dS
+e^{-\tau T}\int_{\Omega\setminus\overline D}
(F_0w-Fw_0)dx
\end{array}
$$
and hence
$$\displaystyle
I_{\partial\Omega}(\tau;B)
=\int_{\partial D}w\frac{\partial w_0}{\partial\nu}dS
+e^{-\tau T}\int_{\Omega\setminus\overline D}
(F_0w-Fw_0)dx.
\tag {2.5}
$$
This is the first representation of the indicator function.  
Next we decompose the first term on the right-hand side of
(2.5).  The result yields the following decomposition formula.

\proclaim{\noindent Proposition 2.1.}
We have
$$\begin{array}{ll}
\displaystyle
I_{\partial\Omega}(\tau;B)
&
\displaystyle
=J(\tau)+E(\tau)+{\cal R}(\tau),
\end{array}
\tag {2.6}
$$
where
$$\displaystyle
J(\tau)=\int_D(\vert\nabla w_0\vert^2+\tau^2\vert w_0\vert^2)dx,
\tag {2.7}
$$
$$\displaystyle
E(\tau)=\int_{\Omega\setminus\overline D}(\vert\nabla R\vert^2+\tau^2\vert R\vert^2)dx
\tag {2.8}
$$
and
$$
\displaystyle
{\cal R}(\tau)
=e^{-\tau T}\left\{\int_DF_0w_0dx
+\int_{\Omega\setminus\overline D}FRdx+\int_{\Omega\setminus\overline D}(F_0-F)w_0dx\right\}.
\tag {2.9}
$$

\endproclaim

{\it\noindent Proof.}
The proof presented here is now standard in the enclosure method, however, in the next section
we make use of an equation appeared in the proof.  So for reader's convenience
we present the proof.

Since $\overline B\cap\overline\Omega=\emptyset$,
the $R$ satisfies
$$\left\{
\begin{array}{ll}
\displaystyle
(\Delta-\tau^2)R=e^{-\tau T}(F-F_0) & \text{in}\,\Omega\setminus\overline D,\\
\\
\displaystyle
\frac{\partial R}{\partial\nu}=0 & \text{on}\,\partial\Omega,\\
\\
\displaystyle
\frac{\partial R}{\partial\nu}=-\frac{\partial w_0}{\partial\nu} & \text{on}\,\partial D.
\end{array}
\right.
\tag {2.10}
$$
Then, one can wite
$$\displaystyle
\int_{\partial D}w\frac{\partial w_0}{\partial\nu}
=\int_{\partial D}w_0\frac{\partial w_0}{\partial\nu}-\int_{\partial D}R\frac{\partial R}{\partial\nu}.
$$
It follows from (2.3) that
$$
\displaystyle
\int_{\partial D}w_0\frac{\partial w_0}{\partial\nu}\,dS
=\int_D(\vert\nabla w_0\vert^2+\tau^2\vert w_0\vert^2)dx
+e^{-\tau T}\int_DF_0w_0dx.
$$
It follows from (2.10) that
$$\begin{array}{ll}
\displaystyle
-\int_{\partial D}R\frac{\partial R}{\partial\nu} dS
&
\displaystyle
=\int_{\partial\,(\Omega\setminus\overline D)}R\frac{\partial R}{\partial\nu} dS\\
\\
\displaystyle
&
\displaystyle
=\int_{\Omega\setminus\overline D}(\vert\nabla R\vert^2+\tau^2\vert R\vert^2)dx
+e^{-\tau T}\int_{\Omega\setminus\overline D}(F-F_0)Rdx.
\end{array}
\tag {2.11}
$$
Thus we obtain
$$\begin{array}{l}
\displaystyle
\,\,\,\,\,\,
\int_{\partial D}w\frac{\partial w_0}{\partial\nu}\\
\\
\displaystyle
=\int_D(\vert\nabla w_0\vert^2+\tau^2\vert w_0\vert^2)dx+
\int_{\Omega\setminus\overline D}(\vert\nabla R\vert^2+\tau^2\vert R\vert^2)dx\\
\\
\displaystyle
\,\,\,
+e^{-\tau T}\left\{\int_DF_0w_0dx+\int_{\Omega\setminus\overline D}(F-F_0)Rdx\right\}.
\end{array}
\tag {2.12}
$$
Then a combination of (2.5) and (2.12) gives (2.6).

\noindent
$\Box$

\subsection{Estimating each term of the decomposition formula}

First we give a {\it rough} estimate of $E(\tau)$ from above in terms of $J(\tau)$.

\proclaim{\noindent Lemma 2.1.}
We have, as $\tau\longrightarrow\infty$
$$\displaystyle
E(\tau)=O\left(\tau^2J(\tau)+\tau^2e^{-2\tau T}\right).
\tag {2.13}
$$
\endproclaim

{\it\noindent Proof.}
It follows from the boundary condition on $\partial D$ in (2.10) and (2.11) that
$$\begin{array}{l}
\displaystyle
\,\,\,\,\,\,
\int_{\Omega\setminus\overline D}\left(\vert\nabla R\vert^2+\tau^2\vert R\vert^2+e^{-\tau T}(F-F_0)R\right)\,dx
\\
\\
\displaystyle
=\int_{\partial D}\frac{\partial w_0}{\partial\nu}RdS,
\end{array}
$$
that is
$$\begin{array}{l}
\displaystyle
\,\,\,\,\,\,
\int_{\Omega\setminus\overline D}
\left(\vert\nabla R\vert^2
+\tau^2\left\vert R+e^{-\tau T}\frac{F-F_0}{2\tau^2}\right\vert^2\right)\,dx\\
\\
\displaystyle
=\int_{\partial D}\frac{\partial w_0}{\partial\nu}RdS
+\frac{e^{-2\tau T}}{4\tau^2}\int_{\Omega\setminus\overline D}\vert F-F_0\vert^2\,dx.
\end{array}
$$
Since from (2.2) and (2.4) we have $\Vert F-F_0\Vert_{L^2(\Omega\setminus\overline D)}=O(\tau)$, this yields
$$\begin{array}{l}
\displaystyle
\,\,\,\,\,\,
E(\tau)
\le 2\int_{\partial D}\frac{\partial w_0}{\partial\nu}RdS+O(e^{-2\tau T}).
\end{array}
\tag {2.14}
$$
By the trace theorem \cite{Gr}, one can choose a positive constant $C=C(D,\Omega)$ and $\tilde{R}\in H^1(D)$
such that $\tilde{R}=R$ on $\partial D$ and $\Vert\tilde{R}\Vert_{H^1(D)}\le C\Vert R\Vert_{H^1(\Omega\setminus\overline D)}$.
Then, we have
$$\begin{array}{l}
\displaystyle
\,\,\,\,\,\,
\int_{\partial D}\frac{\partial w_0}{\partial\nu}RdS\\
\\
\displaystyle
=\int_{\partial D}\frac{\partial w_0}{\partial\nu}\tilde{R}dS\\
\\
\displaystyle
=\int_D(\Delta w_0)\tilde{R}\,dx+\int_D\nabla w_0\cdot\nabla\tilde{R}\,dx\\
\\
\displaystyle
=\tau^2\int_Dw_0\tilde{R}dx+\int_D\nabla w_0\cdot\nabla\tilde{R}\,dx+e^{-\tau T}\int_DF_0\tilde{R}dx.
\end{array}
$$
Note that in the last step, we have made use of equation (2.3) on $D$.
Then the choice of $\tilde{R}$ and (2.4) yield
$$\begin{array}{l}
\displaystyle
\,\,\,\,\,\,
\left\vert\int_{\partial D}\frac{\partial w_0}{\partial\nu}\,R\,dS
\right\vert
\le
C\Vert R\Vert_{H^1(\Omega\setminus\overline D)}
\left(\tau^2\Vert w_0\Vert_{L^2(D)}+\Vert\nabla w_0\Vert_{L^2(D)}+e^{-\tau T}\tau\right).
\end{array}
\tag {2.15}
$$
Here we note that $\Vert R\Vert_{H^1(\Omega\setminus\overline D)}\le E(\tau)^{1/2}$ for all $\tau\ge 1$,
$\Vert w_0\Vert_{L^2(D)}\le\tau^{-1}J(\tau)^{1/2}$, $\Vert\nabla w_0\Vert_{L^2(D)}\le J(\tau)^{1/2}$ for all $\tau>0$.
From these, (2.14) and (2.15) we obtain
$$\begin{array}{l}
\displaystyle
E(\tau)
\le
C'\tau E(\tau)^{1/2}
(J(\tau)^{1/2}+e^{-\tau T})+O(e^{-2\tau T}),
\end{array}
$$
where $C'$ is a positive constant.
Now a standard argument yields (2.13).

\noindent
$\Box$

{\bf\noindent Remark 2.1.}
The advantage of the proof of  Lemma 2.1 is  shown in the right-hand side on (2.15).
We make use of only $H^1$-regularity of $w_0$ in $D$  together with $\Delta w_0\in L^2(D)$.
We do not make use of a concrete expression of the solution of (2.3) at this stage.

Next we describe upper and lower estimates for $J(\tau)$.

\proclaim{\noindent Lemma 2.2.}

(i)  We have, as $\tau\longrightarrow\infty$
$$
\displaystyle
J(\tau)
\displaystyle
=O(\tau^2 e^{-2\tau\text{dist}\,(D,B)}+e^{-2\tau T}).
\tag {2.16}
$$

(ii) Let $T$ satisfies 
$$\displaystyle
T>\text{dist}\,(D,B).
\tag {2.17}
$$
Then, then there exist positive constants $\tau_0$ and $C$ such that, for all $\tau\ge\tau_0$
$$\displaystyle
\tau^{10} e^{2\tau\text{dist}\,(D,B)}J(\tau)\ge C.
\tag {2.18}
$$

\endproclaim

{\it\noindent Proof.}
Set
$$
\displaystyle
\epsilon_0=e^{\tau T}(w_0-v_0),
$$
where $v_0\in H^1(\Bbb R^3)$ is the solution of (1.11).

We have
$$\displaystyle
w_0=v_0+e^{-\tau T}\epsilon_0
$$
and from (2.3)
$$\begin{array}{ll}
\displaystyle
(\Delta-\tau^2)\epsilon_0=F_0 & \text{in}\,\Bbb R^3.
\end{array}
\tag {2.19}
$$
Then, from (2.4) and (2.19) we can easily see that 
$$
\displaystyle \tau\Vert\epsilon_0\Vert_{L^2(\Bbb R^3)}+\Vert\nabla\epsilon_0\Vert_{L^2(\Bbb R^3)}=O(1).
\tag {2.20}
$$
Let $U$ be an arbitrary bounded open subset of $\Bbb R^3$ such that $\overline B\cap\overline U=\emptyset$.
The expression (1.13) for $v_0$ yields
$$
\displaystyle
\tau\Vert v_0\Vert_{L^2(U)}
+
\Vert\nabla v_0\Vert_{L^2(U)}
=O(\tau e^{-\tau\text{dist}\,(U,B)}).
\tag {2.21}
$$
These together with (2.20) give
$$
\displaystyle
\tau\Vert w_0\Vert_{L^2(U)}+
\Vert\nabla w_0\Vert_{L^2(U)}
=O(\tau e^{-\tau\text{dist}\,(U,B)}+e^{-\tau T}).
\tag {2.22}
$$
Now this for $U=D$ and (2.7) yield (2.16).

It follows from (2.20) that
$$
\displaystyle
J(\tau)
\ge \frac{1}{2}J_0(\tau)+O(e^{-2\tau T}),
$$
where
$$\displaystyle
J_0(\tau)=\int_D(\vert\nabla v_0\vert^2+\tau^2\vert v_0\vert^2)dx.
$$
By Lemma A.1 in Appendix we know
$$\displaystyle
J_0(\tau)
\ge
C^2\tau^{-4}\int_D\frac{e^{-2\tau(\vert x-p\vert-\eta)}}{\vert x-p\vert^2}\,dx.
$$
In \cite{IW00, IEE} we have already known that, there exist positive constants $\tau_0$ and $C'$ such that,
for all $\tau\ge\tau_0$ 
$$\displaystyle
\tau^6 e^{2\tau\text{dist}\,(D,B)}\int_D\frac{e^{-2\tau(\vert x-p\vert-\eta)}}{\vert x-p\vert^2}\,dx\ge C'.
$$

Now it is clear that (2.18) is valid under condition (2.17).

\noindent
$\Box$

{\bf\noindent Remark 2.2.}
In the proof of (2.16) the estimate (2.21) is essential.  For the purpose, we made use of 
the explicit expression of $v_0$ given by (1.13).

Now we are ready to give upper bounds for $E(\tau)$ and ${\cal R}(\tau)$.
From (2.13) and (2.16) we obtain
$$
\displaystyle
E(\tau)
=O(\tau^4 e^{-2\tau\text{dist}\,(D,B)}+\tau^2 e^{-2\tau T}).
\tag {2.23}
$$
This yields
$$\displaystyle
\Vert R\Vert_{L^2(\Omega\setminus\overline D)}
=O(\tau e^{-\tau\text{dist}\,(D,B)}+e^{-\tau T}).
\tag {2.24}
$$
This together with (2.2) gives
$$\displaystyle
\int_{\Omega\setminus\overline D}FRdx
=O(\tau^2 e^{-\tau\text{dist}\,(D,B)}+\tau e^{-\tau T}).
\tag {2.25}
$$
And also it follows from (2.2), (2.4) and (2.22) with $U=D, \Omega\setminus\overline D$ we obtain
$$
\displaystyle
\int_DF_0w_0dx=O(\tau e^{-\tau\text{dist}\,(D,B)}+e^{-\tau T})
\tag {2.26}
$$
and
$$
\displaystyle
\int_{\Omega\setminus\overline D}(F_0-F)w_0dx
=O(\tau e^{-\tau\text{dist}\,(\Omega,B)}+e^{-\tau T}).
\tag {2.27}
$$
Applying these to the right-hand side on (2.9), we obtain
$$
\displaystyle
{\cal R}(\tau)
=O(e^{-\tau T}(\tau^2 e^{-\tau\text{dist}\,(D,B)}+\tau e^{-\tau\text{dist}\,(\Omega,B)}
+\tau e^{-\tau T})).
\tag {2.28}
$$

{\bf\noindent Remark 2.3.}
Since $\text{dist}(\Omega,B)<\text{dist}\,(D,B)$,
the estimates (2.25), (2.26) and (2.27) suggest us that the decaying order of the integral in the third term on (2.9) 
is slower than other two terms.  This is a reason from a techinical point of view why we should impose the condition (1.8).

\subsection{Proof of (1.9)}

It follows from (2.28) that
$$\begin{array}{l}
\displaystyle
\,\,\,\,\,\,\,
e^{2\tau\text{dist}\,(D,B)}{\cal R}(\tau)
\\
\\
\displaystyle
=
O(\tau^2 e^{-\tau(T-\text{dist}\,(D,B))}
+\tau e^{-\tau(T-2\text{dist}\,(D,B)+\text{dist}\,(\Omega,B))}
+\tau e^{-2\tau(T-\text{dist}\,(D,B))}).
\end{array}
\tag {2.29}
$$
Now let $T$ satisfy (1.8).
We note that since $\text{dist}\,(D,B)>\text{dist}\,(\Omega,B)$, we have
$$\displaystyle
2\text{dist}\,(D,B)-\text{dist}\,(\Omega,B)>\text{dist}\,(D,B)
$$
and hence
$$\displaystyle
T-(2\text{dist}\,(D,B)-\text{dist}\,(\Omega,B))<T-\text{dist}\,(D,B).
$$
Thus (2.17) is also satisfied.
Then from (2.29) we have, as $\tau\longrightarrow\infty$
$$
\displaystyle
e^{2\tau\text{dist}\,(D,B)}{\cal R}(\tau)=O(\tau^2e^{-c\tau}),
\tag {2.30}
$$
where
$$
c=T-(2\text{dist}\,(D,B)-\text{dist}\,(\Omega,B)).
$$

Rewrite (2.23) as 
$$\displaystyle
e^{2\tau\text{dist}(D,B)}E(\tau)
=O(\tau^4+\tau^2 e^{-2\tau(T-\text{dist}\,(D,B))}).
$$
This gives
$$\displaystyle
e^{2\tau\text{dist}(D,B)}E(\tau)=O(\tau^4).
\tag {2.31}
$$
Moreover, it follows from (2.16) that
$$\displaystyle
e^{2\tau\text{dist}(D,B)}J(\tau)=O(\tau^2).
\tag {2.32}
$$

Now applying (2.30), (2.31) and (2.32) to the right-hand side on (2.6), we obtain
$$\displaystyle
e^{2\tau\text{dist}(D,B)}I_{\partial\Omega}(\tau;B)=O(\tau^4).
\tag {2.33}
$$

Since $E(\tau)\ge 0$, it follows from (2.6) and (2.30) that
$$\displaystyle
\tau^{10}e^{2\tau\text{dist}\,(D,B)}I_{\partial\Omega}(\tau;B)\ge 
\tau^{10} e^{2\tau\text{dist}\,(D,B)}J(\tau)+O(\tau^{12} e^{-c\tau}).
\tag {2.34}
$$
Since $T$ satisfies (2.17), a combination of (2.18) and (2.34) ensures that
there exist positive constants $\tau_0'$ and $C'$ such that, for all $\tau\ge\tau_0$
$$\displaystyle
\tau^{10} e^{2\tau\text{dist}\,(D,B)}I_{\partial\Omega}(\tau;B)\ge C'.
\tag {2.35}
$$
In particular, from this we know that $I_{\partial D}(\tau;B)>0$ for all $\tau\ge\tau_0$ with
a  sufficiently large $\tau_0$.
Now a combination of (2.33) and (2.35) yields (1.9).

\subsection{Proof of (1.10)}

Let $T>0$.
From (2.16), (2.23) and (2.28) we have
$$\left\{
\begin{array}{l}
\displaystyle
e^{\tau T}J(\tau)=O(\tau^2 e^{\tau(T-2\text{dist}\,(D,B))}+e^{-\tau T}),
\\
\\
\displaystyle
e^{\tau T}E(\tau)=O(\tau^4 e^{\tau(T-2\text{dist}\,(D,B))}+\tau^2e^{-\tau T}),
\\
\\
\displaystyle
e^{\tau T}{\cal R}(\tau)
=O(\tau^2 e^{-\tau\text{dist}\,(D,B)}+\tau e^{-\tau\text{dist}\,(\Omega,B)}+\tau e^{-\tau T}).
\end{array}
\right.
\tag {2.36}
$$
Note that there is no restriction on the size of $T$.  Therefore from (2.6) and (2.36) we obtain 
$e^{\tau T}I_{\partial\Omega}(\tau;B)=0$ if $T<2\text{dist}\,(D,B)$.
If $T>2\text{dist}\,(D,B)$, then $T$ satisfies (1.8).  Thus, for sufficiently large $\tau$ we can write
$$\displaystyle
e^{\tau T}I_{\partial\Omega}(\tau;B)
=\exp\,\tau\left(T-2\text{dist}\,(D,B)+\frac{1}{\tau}\log I_{\partial\Omega}(\tau;B)+2\text{dist}\,(D,B)\right).
$$
Then it follows from (1.9) that $\lim_{\tau\longrightarrow\infty}e^{\tau T}I_{\partial\Omega}(\tau;B)=\infty$.

This completes the proof of (1.10).

\subsection{Proof of Corollaries}

The proof of Corollary 1.1 is as follows.

It follows from (2.23) that
$$\displaystyle
\Vert w-w_0\Vert_{H^{1/2}(\partial\Omega)}=O(\tau^2 e^{-\tau\text{dist}\,(D,B)}
+\tau e^{-\tau T}
).
\tag {2.37}
$$
A combination of the standard estimate
$$\displaystyle
\left\Vert\frac{\partial\epsilon_0}{\partial\nu}\right\Vert_{H^{-1/2}(\partial\Omega)}
\le C\left(\Vert\Delta\epsilon_0\Vert_{L^2(\Omega)}+\Vert\nabla\epsilon_0\Vert_{L^2(\Omega)}\right)
$$
together with (2.4), (2.19) and (2.20), we obtain
$$
\displaystyle
\left\Vert\frac{\partial\epsilon_0}{\partial\nu}\right\Vert_{H^{-1/2}(\partial\Omega)}=O(\tau).
$$
A combination of this and (2.37) gives
$$\displaystyle
\left\vert\int_{\partial\Omega}(w-w_0)\frac{\partial\epsilon_0}{\partial\nu}dS\right\vert
=O(\tau^3e^{-\tau\text{dist}\,(D,B)}
+\tau^2e^{-\tau T}
).
$$
Hence we obtain
$$\displaystyle
\displaystyle
I_{\partial\Omega}^s(\tau;B)=I_{\partial\Omega}(\tau;B)
+O(\tau^3e^{-\tau(T+\text{dist}\,(D,B))}+\tau^2e^{-2\tau T}
).
$$
Then, we can easily check all the statements of Theorem 1.1 are transplanted into those of Corollary 1.1
by using the following facts.

$\bullet$  One can write
$$\left\{\begin{array}{l}
\displaystyle
e^{2\tau\text{dist}\,(D,B)}\tau^3e^{-\tau(T+\text{dist}\,(D,B))}
=\tau^3e^{-\tau(T-\text{dist}\,(D,B))},\\
\\
\displaystyle
e^{2\tau\text{dist}\,(D,B)}\tau^2e^{-2\tau T}
=\tau^2 e^{-2\tau(T-\text{dist}\,(D,B))}.
\end{array}
\right.
$$

$\bullet$  (1.8) implies
$$\displaystyle
T-\text{dist}\,(D,B)>\text{dist}\,(D,B)-\text{dist}\,(\Omega,B)
$$
and we have $\text{dist}\,(D,B)>\text{dist}\,(\Omega,B)$.

$\bullet$  One has
$$\left\{\begin{array}{l}
\displaystyle
e^{\tau T}\tau^3e^{-\tau(T+\text{dist}\,(D,B))}
=\tau^3 e^{-\tau\text{dist}\,(D,B)},\\
\\
\displaystyle
e^{\tau T}\tau^2 e^{-2\tau T}
=\tau^2 e^{-\tau T}.
\end{array}
\right.
$$

The proof of Corollary 1.2 is as follows.
From the expression (1.13), we have
$$\displaystyle
\left\Vert\frac{\partial v_0}{\partial\nu}\right\Vert_{L^2(\partial\Omega\setminus\partial\Omega(M,B))}
=O(\tau e^{-\tau M}).
\tag {2.38}
$$
It follows from (2.37) that
$$\displaystyle
\Vert w-w_0\Vert_{L^2(\partial\Omega)}
=O(\tau^2 e^{-\tau\text{dist}\,(D,B)}+\tau e^{-\tau T}
).
$$
A combination of these gives
$$\displaystyle
I_{\partial\Omega}(\tau;B,M)=I_{\partial\Omega}^s(\tau;B)+O(\tau^3e^{-\tau\text{dist}\,(D,B)}e^{-\tau M}
+\tau^2 e^{-\tau T}e^{-\tau M}
).
$$
It is clear that we can check the validity of the statement in Corollary 1.2
by using the following facts.

$\bullet$  One can write
$$\left\{\begin{array}{l}
\displaystyle
e^{2\tau\text{dist}\,(D,B)}\tau^3e^{-\tau\text{dist}\,(D,B)}e^{-\tau M}
=\tau^3 e^{-\tau(M-\text{dist}\,(D,B))},\\
\\
\displaystyle
e^{2\tau\text{dist}\,(D,B)}\tau^2 e^{-\tau T}e^{-\tau M}
=\tau^2 e^{-\tau(T+M-2\text{dist}\,(D,B))}.
\end{array}
\right.
$$

$\bullet$  (1.15) implies
$$\begin{array}{l}
\displaystyle
\,\,\,\,\,\,T+M-2\text{dist}\,(D,B)
\\
\\
\displaystyle
\ge
2M-\text{dist}\,(\Omega,B)+M-2\text{dist}\,(D,B)
\\
\\
\displaystyle
=2(M-\text{dist}\,(D,B))+(M-\text{dist}\,(\Omega,B))
\\
\\
\displaystyle
>2(M-\text{dist}\,(D,B))+(M-\text{dist}\,(D,B))\\
\\
\displaystyle
=3(M-\text{dist}\,(D,B)).
\end{array}
$$

{\bf\noindent Remark 2.4.}  In the proof of Corollary 1.2 the estimate (2.38) is essential.
This also employs the explicit expression (1.13) for $v_0$.

\section{Conclusion}

We have {\it indicated} the idea of making use of a special, however, natural Neumann data 
like (1.4) by using a proto-type inverse obstacle problem. 
The Neumann data tell us: {\it how to hit} the surface of a body to obtain the information
about the distance of an arbitrary given point outside the body to an unknown obstacle.
The Neumann data are given by solving the initial value problem for the governing equation
with special initial data in the whole space and taking the Neumann derivative
on the surface of the body.  For constant coefficient case it will be possible to
obtain the data explicitly.

In principle, it would be possible to apply the idea developed in this paper to other
inverse obstacle problems using time domain data over a finite time interval whose
governing equations are given by hyperbolic or parabolic equations/systems in a {\it bounded/unbounded spacial domain}
with an {\it outer boundary}, including half-space, slab, etc. 
The program of realizing the idea shall be done step by step in forthcoming papers.

$$\quad$$

\centerline{{\bf Acknowledgment}}

The author was partially supported by Grant-in-Aid for
Scientific Research (C)(No. 25400155) of Japan  Society for
the Promotion of Science.

$$\quad$$

\section{Appendix}

In this appendix we give a proof of the following Lemma.

\proclaim{\noindent Lemma A.1.}
There exist positive constants $C$ and $\tau_0$ such that, for all $\tau\ge\tau_0$ and all $x\in\Bbb R^3\setminus\overline B$
$$\displaystyle
v_0(x)\ge C\tau^{-3}\frac{e^{-\tau(\vert x-p\vert-\eta)}}{\vert x-p\vert}.
\tag {A.1}
$$

\endproclaim

The proof presented here is based on an explicit computation of a volume integral.

\proclaim{\noindent Proposition A.1.}
Let $x\in\Bbb R^3\setminus\overline B$.
We have
$$
\displaystyle
\int_B\frac{e^{-\tau\vert x-y\vert}}{\vert x-y\vert}\cdot\vert y-p\vert\,dy
=\frac{4\pi}{\tau^2}
\frac{e^{-\tau\vert x-p\vert}}{\vert x-p\vert}\left\{\left(\eta+\frac{2}{\tau^2}\right)\cosh\tau\eta
-\frac{2\eta}{\tau}\sinh\tau\eta
-\frac{2}{\tau^2}\right\}.
\tag {A.2}
$$

\endproclaim

{\it\noindent Proof.}  We denote by $I(x)$ the left-hand side on (A.2).
It suffices to consider the case when $p=0$.
The change of variables $y=r\omega\,(0<r<\eta,\,\omega\in S^2)$ and a rotation give us
$$\begin{array}{ll}
\displaystyle
I(x) & \displaystyle =\int_0^{\eta}r^3 dr\int_{S^2}\frac{e^{-\tau\vert x-r\omega\vert}}{\vert x-r\omega\vert}d\omega
\\
\\
\displaystyle
& \displaystyle =\int_0^{\eta}r^3 dr\int_{S^2}
\frac{\displaystyle e^{-\tau\vert \vert x\vert\mbox{\boldmath $e$}_3-r\omega\vert}}
{\displaystyle\vert\vert x\vert\mbox{\boldmath $e$}_3-r\omega\vert}d\omega\\
\\
\displaystyle
& \displaystyle =\int_0^{\eta}r^3dr\int_0^{2\pi}d\theta
\int_0^{\pi}\sin\varphi d\varphi
\frac{\displaystyle e^{-\tau\sqrt{\vert x\vert^2-2r\vert x\vert\cos\varphi+r^2}}}
{\displaystyle\sqrt{\vert x\vert^2-2r\vert x\vert\cos\varphi+r^2}}\\
\\
\displaystyle
& \displaystyle =2\pi\int_0^{\eta}Q(\vert x\vert,r)r^3 dr,
\end{array}
\tag {A.3}
$$
where $\mbox{\boldmath $e$}_3=(0,0,1)$ and 
$$\begin{array}{lll}
\displaystyle
Q(\xi,r)
=\int_0^{\pi}
\frac{\displaystyle e^{-\tau\sqrt{\xi^2-2r\xi\cos\varphi+r^2}}}
{\displaystyle\sqrt{\xi^2-2r\xi\cos\varphi+r^2}}\sin\varphi d\varphi, & \xi>\eta, & 0<r<\eta.
\end{array}
$$
Fix $\xi\in]\eta,\,\infty [$ and $r\in]0,\,\eta[$.
The change of variable
$$\displaystyle
s=\sqrt{\xi^2-2r\xi\cos\varphi+r^2},\,\varphi\in]0,\,\pi[
$$
gives
$$\displaystyle
s^2=\xi^2-2r\xi\cos\varphi+r^2
$$
and
$$\displaystyle
sds=r\xi\sin\varphi d\varphi.
$$
Thus we have
$$\begin{array}{ll}
\displaystyle
Q(\xi,r)
& \displaystyle
=\frac{1}{r\xi}\int_{\xi-r}^{\xi+r}e^{-\tau s}ds
\\
\\
\displaystyle
& \displaystyle =-\frac{1}{r\xi\tau}\left(e^{-\tau(\xi+r)}-e^{-\tau(\xi-r)}\right).
\end{array}
$$
From this and (A.3) we obtain
$$
\displaystyle
I(x)=\frac{2\pi}{\xi\tau}\int_0^{\eta}\left(e^{-\tau(\xi-r)}-e^{-\tau(\xi+r)}\right)r^2dr\vert_{\xi=\vert x\vert}.
\tag {A.4}
$$
Since we have
$$\displaystyle
\left\{
\begin{array}{l}
\displaystyle
\int_0^{\eta}e^{-\tau(\xi-r)}r^2dr
=\frac{1}{\tau}
\left(\eta^2-\frac{2\eta}{\tau}+\frac{2}{\tau^2}\right)e^{-\tau(\xi-\eta)}
-\frac{2}{\tau^3}e^{-\tau\xi},
\\
\\
\displaystyle
\int_0^{\eta}e^{-\tau(\xi+r)}r^2dr
=-\frac{1}{\tau}
\left(\eta^2+\frac{2\eta}{\tau}+\frac{2}{\tau^2}\right)e^{-\tau(\xi+\eta)}
+\frac{2}{\tau^3}e^{-\tau\xi},
\end{array}
\right.
$$
one gets
$$
\displaystyle
\int_0^{\eta}\left(e^{-\tau(\xi-r)}-e^{-\tau(\xi+r)}\right)r^2dr
=\frac{2}{\tau}e^{-\tau\xi}
\left\{\left(\eta+\frac{2}{\tau^2}\right)\cosh\tau\eta
-\frac{2\eta}{\tau}\sinh\tau\eta
-\frac{2}{\tau^2}\right\}.
$$
Now from this and (A.4) we obtain (A.2).

\noindent
$\Box$

Write
$$\begin{array}{l}
\displaystyle
I(x)
=\frac{4\pi}{\tau^4}
\frac{e^{-\tau\vert x-p\vert}}{\vert x-p\vert}
\left\{\left((\tau\eta)^2+2\right)\cosh\tau\eta
-2\tau\eta\sinh\tau\eta
-2\right\}.
\end{array}
$$
Define
$$\displaystyle
\Psi(s)
=\left(s^2+2\right)\cosh s
-2s\sinh s
-2.
$$
Then we have the expression
$$
\displaystyle
I(x)=\frac{4\pi}{\tau^4}\frac{e^{-\tau\vert x-p\vert}}{\vert x-p\vert}
\Psi(\tau\eta).
\tag {A.5}
$$

We know that 
$$\begin{array}{ll}
\displaystyle
\int_B\frac{e^{-\tau\vert x-y\vert}}{\vert x-y\vert}\,dy=\frac{4\pi\varphi(\tau\eta)}{\tau^3}\frac{e^{-\tau\vert x-p\vert}}{\vert x-p\vert}, & x\in\Bbb R^3\setminus\overline B,
\end{array}
\tag {A.6}
$$
where
$$\displaystyle
\varphi(s)=s\cosh s-\sinh s.
$$
This is nothing but a consequence of the mean value theorem for the modified Helmholtz equation
\cite{CH} and can be checked directly by using a similar computation as above.

Thus (A.5) and (A.6) give us the expression for $v_0(x)$ for $x\in\Bbb R^3\setminus\overline B$:
$$\displaystyle
v_0(x)
=\frac{1}{\tau^4}
\frac{e^{-\tau\vert x-p\vert}}{\vert x-p\vert}
\left(s\varphi(s)-\Psi(s)\right)\vert_{s=\tau\eta}.
$$
Since we have, as $s\longrightarrow\infty$
$$\displaystyle
s\varphi(s)-\Psi(s)
=-2\cosh s+s\sinh s+2\sim se^{s},
$$
one gets (A.1).

\vskip1cm
\noindent
e-mail address

ikehata@amath.hiroshima-u.ac.jp

\end{document}